\documentclass[11pt,a4paper]{amsart}

\usepackage{amsmath,amssymb,amsfonts,amsthm,mathtools}
\usepackage{enumitem}
\usepackage{booktabs,array}
\usepackage[hidelinks]{hyperref}
\hypersetup{
  pdftitle={Cyclic crossed-product envelopes and arity-detecting relative K-theory of automorphism-derived n-ary algebras},
  pdfauthor={Chandrasekhar Gokavarapu, Sajani Lavanya Madasi, Madhusudhana Rao Dasari}
}

\newtheorem{theorem}{Theorem}[section]
\newtheorem{proposition}[theorem]{Proposition}
\newtheorem{lemma}[theorem]{Lemma}
\newtheorem{corollary}[theorem]{Corollary}
\theoremstyle{definition}
\newtheorem{definition}[theorem]{Definition}
\newtheorem{example}[theorem]{Example}
\theoremstyle{remark}
\newtheorem{remark}[theorem]{Remark}

\DeclareMathOperator{\Aut}{Aut}
\DeclareMathOperator{\Ad}{Ad}
\DeclareMathOperator{\End}{End}
\DeclareMathOperator{\hofib}{hofib}
\DeclareMathOperator{\cofib}{cofib}
\DeclareMathOperator{\Stab}{Stab}
\DeclareMathOperator{\Vect}{Vect}

\newcommand{\K}{\mathbf K}

\newcommand{\Ua}{\mathsf U^{\mathrm{an}}}
\newcommand{\Uac}{\mathsf U^{\mathrm{an}}_c}
\newcommand{\R}{\mathbf R}

\newcommand{\op}{\mathrm{op}}
\newcommand{\id}{\mathrm{id}}

\newcommand{\eps}{\varepsilon}
\newcommand{\doilink}[1]{\href{https://doi.org/#1}{\nolinkurl{doi:#1}}}
\newcommand{\mrlink}[1]{\href{https://www.ams.org/mathscinet-getitem?mr=#1}{MathSciNet MR#1}}

\title[Cyclic crossed-product envelopes]{Cyclic crossed-product envelopes and arity-detecting relative $K$-theory of automorphism-derived $n$-ary algebras}

\author{Chandrasekhar Gokavarapu}
\address{Department of Mathematics, Acharya Nagarjuna University, Guntur, Andhra Pradesh 522510, India}
\address{Department of Mathematics, Government College (Autonomous), Rajahmundry, Andhra Pradesh 533105, India}

\email{chandrasekhargokavarapu@gmail.com}
\email{chandrasekhargokavarapu@gcrjy.ac.in}
\urladdr{https://orcid.org/0009-0006-5306-371X}
\thanks{Corresponding author: Chandrasekhar Gokavarapu.}

\author{Sajani Lavanya Madasi}
\address{Department of Mathematics, Government College (Autonomous), Rajahmundry, Andhra Pradesh 533105, India}
\email{sajanimaths@gcrjy.ac.in}
\urladdr{https://orcid.org/0009-0000-8194-4027}

\author{Madhusudhana Rao Dasari}
\address{Department of Mathematics, Acharya Nagarjuna University, Guntur, Andhra Pradesh 522510, India}
\address{Department of Mathematics, Government College for Women (Autonomous), Guntur, Andhra Pradesh 522006, India}

\email{dmrmaths@gmail.com}
\urladdr{https://orcid.org/0000-0003-0897-868X}

\date{}

\begin{document}

\begin{abstract}
Let $m=n-1$.  For a unital associative algebra $A$, an automorphism $\sigma$,
and $c\in A^\times$ satisfying $\sigma(c)=c$ and
$\sigma^m=\operatorname{Ad}_c$, we study the Hossz\'u--Gluskin-type product
$[a_0,\ldots,a_m]=a_0\sigma(a_1)\cdots\sigma^{m-1}(a_{m-1})ca_m$.
We prove that its endpoint-anchored all-slot envelope is the cyclic crossed
product generated by $A\otimes A^{\mathrm{op}}$ and $u$, with
$ux=\sigma^e(x)u$ and $u^m=c\otimes(c^{-1})^{\mathrm{op}}$.
For the finite-order case $c=1$ and permutation algebras $A=k^X$, this gives
an orbit--stabilizer block decomposition.  A general dimension-defect theorem
then computes the homotopy fibre of the forgetful $K$-theory map.  In
characteristic zero, transitive $C_m$-sets of every divisor length $d\mid m$
have the same ordinary spectrum $\K(k)^m$ but boundary $\mathbb Z/d$.
The regular action gives $\pi_{-1}\cong\mathbb Z/m$, recovering the arity;
over finite fields, the remaining torsion is governed by multiplicative
orders.
\end{abstract}

\maketitle

\section{Introduction}

A central difficulty in the representation theory of an $n$-ary algebra is
that a module variable may occupy any one of the $n$ input positions.  The
operators generated by these different positions need not collapse to the
left and right multipliers of an ordinary associative algebra.  This paper
constructs a family in which the positional information survives as a finite
cyclic symmetry and can be detected by algebraic $K$-theory.

Totally associative $n$-ary algebras and their modules have been studied from
operadic and cohomological viewpoints
\cite{Gnedbaye1995,Gnedbaye1997,MichorVinogradov1996,Balavoine1998,
AtaguemaMakhlouf2009,BagherzadehBremner2020}, as well as through universal
covers \cite{Sitarz1998}.  Automorphism-derived polyadic products are classical
in the Hossz\'u--Gluskin tradition \cite{DudekGlazek2008}.  Representation
theories for polyadic groups have also been developed
\cite{BorowiecDudekDuplij2006,DudekShahryari2012}.
General operadic theory assigns an associative enveloping algebra whose
modules represent modules over an algebra governed by an operad
\cite{Fresse2009,LodayVallette2012,Khoroshkin2022}.  Polyadic analogues of the
Grothendieck completion have also been proposed \cite{Duplij2022}, and a formal
higher-$K$-theory programme for slot-sensitive $n$-ary $\Gamma$-semiring
modules was developed in \cite{Gokavarapu2025Higher,Gokavarapu2025Fundamental}.
Those constructions concern, respectively, polyadic completion groups and
exact/Waldhausen categories attached to different algebraic objects; they do
not identify the anchored envelope or compute the relative spectrum of the
forgetful functor considered here.  Accordingly, the product below is not
claimed as a new polyadic construction.  The new problem addressed here is to
determine its all-slot module envelope explicitly and to decide whether exact
$K$-theory of the forgetful functor retains the arity.

Let $m=n-1$.  We begin with the full inner-periodic datum: a unital
associative $k$-algebra $A$, an automorphism $\sigma$, and an invertible element
$c$ such that
\[
 \sigma(c)=c,
 \qquad
 \sigma^m=\Ad_c.
\]
The associated product is
\begin{equation}\label{eq:intro-product}
 [a_0,a_1,\ldots,a_m]_{\sigma,c}
 =a_0\,\sigma(a_1)\cdots\sigma^{m-1}(a_{m-1})c a_m.
\end{equation}
This is the associative-algebra analogue of the full Hossz\'u--Gluskin
construction, rather than only its finite-order specialization.  We prove
that every evaluation of a word of length $1+qm$ is the same ordered twisted
product with a terminal factor $c^q$.

A slot module is called \emph{endpoint-anchored} when the two contexts obtained
by placing $c^{-1}$ at the endpoint opposite the module variable act as the
identity.  The first principal result identifies its representing algebra.
Put $A^e=A\otimes_kA^{\op}$, $\sigma^e=\sigma\otimes\sigma^{\op}$, and
$z_c=c\otimes(c^{-1})^{\op}$.  Then
\begin{equation}\label{eq:intro-envelope}
 \Uac(A,\sigma,n)
 \cong
 A^e\langle u\rangle/
 \bigl(ux-\sigma^e(x)u,\ u^m-z_c\bigr).
\end{equation}
The compatibility conditions imply $(\sigma^e)^m=\Ad_{z_c}$ and
$\sigma^e(z_c)=z_c$, so the right-hand side is a cyclic crossed product.  The
proof gives a context normal form: if the hole of a flattened context occurs
at position $h=tm+r$, then the context reduces to a left multiplier, a right
multiplier, and $u^r$, with the carry $c^t$ absorbed into the two multipliers.
Consequently endpoint-anchored slot modules are precisely $A$-bimodules $M$
with an invertible semilinear operator $S$ satisfying
\[
 S(axb)=\sigma(a)S(x)\sigma(b),
 \qquad
 S^m(x)=cxc^{-1}.
\]
For $c=1$ this becomes the skew group algebra
$(A\otimes A^{\op})\rtimes C_m$.  We write $T(A,\sigma)=T(A,\sigma,1)$
in this finite-order case.  Thus a genuinely higher-arity
representation problem becomes cyclically equivariant bimodule theory, while
the inner-periodic theorem also covers automorphisms whose $m$-th power is
nontrivial but inner.

The second principal result treats permutation algebras.  If a cyclic group
$C_m$ acts on a finite set $X$ and $A=k^X$, then the envelope is the algebra of
the action groupoid of the diagonal action on $X\times X$.  We prove the block
formula
\begin{equation}\label{eq:intro-orbit}
 \Ua(k^X,\sigma,n)
 \cong
 \prod_{O\in (X\times X)/C_m}
 M_{|O|}\bigl(k\Stab(x_O)\bigr).
\end{equation}
Write the $C_m$-orbits of $X$ as $C_m/H_\alpha$ and put
$\ell_\alpha=[C_m:H_\alpha]$.  A direct orbit calculation shows that the
pair $(\alpha,\beta)$ contributes $\gcd(\ell_\alpha,\ell_\beta)$ diagonal
orbits, each of length $\operatorname{lcm}(\ell_\alpha,\ell_\beta)$ and with
stabilizer $H_\alpha\cap H_\beta$.  In particular,
\begin{equation}\label{eq:intro-dx}
 d_X:=\gcd\{|O|:O\in(X\times X)/C_m\}
 =\gcd_\alpha\ell_\alpha
 =[C_m:\langle H_\alpha:\alpha\in I\rangle].
\end{equation}
Thus the boundary invariant can be read from the original action rather than
from a separate diagonal-orbit enumeration.

The third principal result is a general dimension-defect theorem.  If a
finite-length $k$-linear abelian category has $N$ absolutely simple objects
and an exact functor to finite-dimensional vector spaces sends them to
spaces of dimensions $w_1,\ldots,w_N$, then, noncanonically,
\begin{equation}\label{eq:intro-dimension-defect}
 \hofib\K(F)\simeq
 \K(k)^{N-1}\times\Omega\K(k;\mathbb Z/\gcd(w_1,\ldots,w_N)).
\end{equation}
Applying this to \eqref{eq:intro-orbit}, and writing
$N_X=\sum_O|\Stab(x_O)|_{p'}$ over an algebraically closed field of
characteristic $p\geq0$, gives
\begin{equation}\label{eq:intro-orbit-fibre}
 \R_X(k)\simeq
 \K(k)^{N_X-1}\times\Omega\K(k;\mathbb Z/d_X).
\end{equation}
This information is not contained in ordinary exact $K$-theory.  In
characteristic $0$, for every divisor $d\mid m$ the transitive action on a set
of size $d$ has ordinary spectrum $\K(k)^m$, while its relative boundary is
$\mathbb Z/d$.  Thus the entire divisor lattice of $m$ is visible to the
forgetful-functor fibre but invisible to ordinary exact $K$-theory.

For the regular action $X=C_m$, the diagonal action has $m$ free orbits, so
\begin{equation}\label{eq:intro-matrix}
 \Ua(k^{C_m},\sigma,m+1)
 \cong \prod_{d\in C_m}M_m(k).
\end{equation}
The corresponding $n$-ary algebra has no polyadic unit when $m>1$.  Hence this
example is not the $n$-fold product of a unital binary algebra, and the cyclic
operator cannot be removed by unit normalization.  Here $N_X=d_X=m$, and the
general orbit invariant becomes arity-sensitive.  For the exact category
$\mathcal M_m$ of anchored modules and the forgetful functor
$F_m:\mathcal M_m\to\Vect_k^{\mathrm{fd}}$, we prove
\begin{equation}\label{eq:intro-fibre}
 \hofib\bigl(\K(F_m)\bigr)
 \simeq
 \K(k)^{m-1}\times\Omega\K(k;\mathbb Z/m),
\end{equation}
where $\K(k;\mathbb Z/m)$ is the cofiber of multiplication by $m$ on
$\K(k)$.  In particular,
\begin{equation}\label{eq:intro-boundary}
 \pi_{-1}\hofib\bigl(\K(F_m)\bigr)\cong\mathbb Z/m.
\end{equation}
Thus the relative spectrum recovers $m=n-1$ and therefore the arity.  This is
the opposite of the familiar unit-derived situation, where a normalized
operator envelope may become independent of $n$.

For a finite field $k=\mathbb F_q$, Quillen's calculation of $K_*(\mathbb
F_q)$ gives all homotopy groups explicitly.  If
$g_j=\gcd(m,q^j-1)$, then
\[
 \pi_{2j-1}\cong
 (\mathbb Z/(q^j-1))^{m-1}\oplus\mathbb Z/g_j,
 \qquad
 \pi_{2j}\cong\mathbb Z/g_{j+1},
\]
with the separate formulas for degrees $0$ and $-1$ stated below.  The cyclic
factor is therefore visible throughout the relative theory, not merely on
$K_0$.

The construction is related to, but distinct from, Sitarz's
$\mathbb Z/(n-1)$-graded universal cover of an $n$-ary algebra
\cite{Sitarz1998}.  The universal cover contains the algebra as a homogeneous
component; the envelope studied here represents an anchored module category.
The appearance of $C_{n-1}$ in both settings reflects the same arity modulus,
but the universal properties and the resulting invariants differ.

The paper is organized as follows.  Section~\ref{sec:twist} proves total associativity for the inner-periodic twist and introduces endpoint-anchored slot modules.  Section~\ref{sec:envelope} proves the cyclic crossed-product envelope theorem.  Section~\ref{sec:permutation} gives
the orbit--stabilizer decomposition and the regular cyclic example.
Section~\ref{sec:K} proves the dimension-defect theorem and computes the relative spectrum, and
Section~\ref{sec:finite-fields} gives the finite-field groups and small-arity
examples.

\section{Inner-periodic automorphism-derived algebras}
\label{sec:twist}

Throughout, $k$ is a commutative unital ring, $m\geq2$, and $n=m+1$.
Let $A$ be a unital associative $k$-algebra, let $c\in A^\times$, and let
$\sigma\in\Aut_k(A)$ satisfy
\begin{equation}\label{eq:inner-periodic-data}
 \sigma(c)=c,
 \qquad
 \sigma^m(a)=cac^{-1}\quad(a\in A).
\end{equation}

\begin{definition}\label{def:twisted}
The \emph{inner-periodic automorphism-derived $n$-ary algebra}
$T(A,\sigma,c)$ is the underlying $k$-module of $A$ with multiplication
\begin{equation}\label{eq:twisted}
 \mu_{\sigma,c}(a_0,\ldots,a_m)
 :=a_0\sigma(a_1)\cdots\sigma^{m-1}(a_{m-1})ca_m.
\end{equation}
Equivalently,
$\mu_{\sigma,c}(a_0,\ldots,a_m)
 =(\prod_{j=0}^{m}\sigma^j(a_j))c$.
\end{definition}

\begin{theorem}[Inner-periodic twist theorem]\label{thm:total-assoc}
The multiplication \eqref{eq:twisted} is totally associative.  More generally,
every parenthesized word of length $1+qm$ has the common value
\begin{equation}\label{eq:long-word}
 \left(\prod_{j=0}^{qm}\sigma^j(a_j)\right)c^q,
\end{equation}
where factors occur in increasing order of $j$.
\end{theorem}

\begin{proof}
It is enough to compare the $m+1$ elementary rebracketings of a word of
length $2m+1$.  Insert the inner product in outer position $i$.  Its image
under $\sigma^i$ is
\[
 \left(\prod_{r=0}^{m}\sigma^{i+r}(a_{i+r})\right)c,
\]
because $\sigma(c)=c$.  For each entry after the inner block, the outer
exponent is smaller than its absolute exponent by $m$.  The identity
$\sigma^m(x)=cxc^{-1}$ moves the intervening $c$ across the remaining factors;
explicitly,
\[
 c\,x_1\cdots x_s\,c
 =(cx_1c^{-1})\cdots(cx_sc^{-1})c^2.
\]
Hence every rebracketing equals
$(\prod_{j=0}^{2m}\sigma^j(a_j))c^2$.  Induction on the number of vertices,
or equivalently successive elementary rebracketings, proves
\eqref{eq:long-word}.
\end{proof}

\begin{definition}\label{def:slot-module}
A \emph{slot module} over $T(A,\sigma,c)$ is a $k$-module $M$ with multilinear
maps
\[
 \lambda_i:A^{\otimes i}\otimes M\otimes A^{\otimes(m-i)}\longrightarrow M,
 \qquad 0\leq i\leq m,
\]
such that all evaluations of a planar $n$-ary word containing exactly one
entry of $M$ agree under total-associativity moves.  It is
\emph{endpoint-anchored} if
\begin{equation}\label{eq:c-anchor}
 \lambda_0(x,1,\ldots,1,c^{-1})=x,
 \qquad
 \lambda_m(c^{-1},1,\ldots,1,x)=x.
\end{equation}
Morphisms commute with every slot operation.
\end{definition}

The regular slot module $M=A$ is endpoint-anchored.  When $c=1$, the two
conditions in \eqref{eq:c-anchor} reduce to the endpoint unit normalizations
used for the finite-order cyclic twist.

\section{The cyclic crossed-product envelope}
\label{sec:envelope}

A one-hole context is a finite planar rooted full $n$-ary tree whose ordinary
leaves are labelled by elements of $A$ and whose distinguished leaf is
labelled by $\star$.  Take the free $k$-module on contexts and impose
multilinearity, scalar contraction by \eqref{eq:twisted}, and
total-associativity.  Grafting makes the quotient an associative algebra.
We further impose the two endpoint anchors \eqref{eq:c-anchor}; the resulting
algebra is denoted $\Uac(A,\sigma,n)$.  Its unital modules are exactly the
endpoint-anchored slot modules.  This is a concrete all-slot realization of
the operadic enveloping construction
\cite{Fresse2009,LodayVallette2012,Khoroshkin2022}.

Put $A^e=A\otimes_kA^{\op}$, $\sigma^e=\sigma\otimes\sigma^{\op}$ and
\[
 z_c=c\otimes(c^{-1})^{\op}\in(A^e)^\times.
\]
Then $(\sigma^e)^m=\Ad_{z_c}$ and $\sigma^e(z_c)=z_c$.  Define the
\emph{cyclic crossed product} as the left $A^e$-module
\begin{equation}\label{eq:cyclic-crossed-product}
 \mathcal C(A,\sigma,c)=\bigoplus_{r=0}^{m-1}A^eu^r
\end{equation}
with multiplication
\begin{equation}\label{eq:crossed-product-multiplication}
 (xu^r)(yu^s)=x(\sigma^e)^r(y)z_c^{\delta}u^{r+s-\delta m},
 \qquad
 \delta=\left\lfloor\frac{r+s}{m}\right\rfloor.
\end{equation}
The identities $(\sigma^e)^m=\Ad_{z_c}$ and $\sigma^e(z_c)=z_c$ make this
multiplication associative.  Equivalently, it has the presentation
$A^e\langle u\rangle/(ux-\sigma^e(x)u,\ u^m-z_c)$.  If $c=1$, it is the
skew group algebra $A^e\rtimes_{\sigma^e}C_m$.

For $a,b\in A$, let $L_a$ be the elementary context with $ac^{-1}$ in
position $0$, the hole in position $m$, and units elsewhere.  Let $R_b$ have
the hole in position $0$, $c^{-1}b$ in position $m$, and units elsewhere.
Let $S$ be the elementary context with the hole in position $1$, $c^{-1}$ in
position $m$, and units elsewhere.  The anchor relations make $L_1=R_1=1$.

\begin{lemma}[Canonical crossed-product relations]\label{lem:crossed-relations}
In $\Uac(A,\sigma,n)$ one has
\begin{equation}\label{eq:crossed-relations}
\begin{gathered}
 L_aL_d=L_{ad},\qquad R_bR_e=R_{eb},\qquad L_aR_b=R_bL_a,\\
 SL_a=L_{\sigma(a)}S,\qquad SR_b=R_{\sigma(b)}S,\qquad
 S^m=L_cR_{c^{-1}}.
\end{gathered}
\end{equation}
\end{lemma}

\begin{proof}
For the first three identities, graft the corresponding elementary contexts,
flatten the resulting two-vertex word, and use the endpoint anchors to remove
the two occurrences of $c^{-1}$; the remaining labels multiply as $ad$ on
the left and as $eb$ on the right.  Moving a left or right multiplier through
the one-step context $S$ shifts its absolute leaf position by one, hence
applies $\sigma$ and gives the two covariance identities.  Finally, flattening
$m$ copies of $S$ places the hole $m$ positions later.  Formula
\eqref{eq:long-word}, together with $\sigma^m=\Ad_c$, changes this displacement
into left multiplication by $c$ and right multiplication by $c^{-1}$.
Thus $S^m=L_cR_{c^{-1}}$ as a context identity.
\end{proof}

\begin{lemma}[Carried context normal form]\label{lem:normal-form}
Let $C$ be a one-hole context with $q$ vertices.  Flatten its leaves as
\[
 x_0,\ldots,x_{h-1},\star,x_{h+1},\ldots,x_{qm},
 \qquad h=tm+r,
 \quad 0\le r<m.
\]
Set
\begin{align}
 a_C&=\left(\prod_{j<h}\sigma^j(x_j)\right)c^t,
 \label{eq:normal-a}\\
 b_C&=c^{-t}\left(\prod_{j>h}\sigma^j(x_j)\right)c^q.
 \label{eq:normal-b}
\end{align}
Then
\begin{equation}\label{eq:normal-form}
 [C]=L_{a_C}R_{b_C}S^r
 \qquad\text{in }\Uac(A,\sigma,n).
\end{equation}
Indeed, the right-hand side has formal action
\begin{equation}\label{eq:normal-form-action}
 a_CS^r(x)b_C
 =\left(\prod_{j<h}\sigma^j(x_j)\right)
  \sigma^h(x)
  \left(\prod_{j>h}\sigma^j(x_j)\right)c^q.
\end{equation}
\end{lemma}

\begin{proof}
For an elementary context, \eqref{eq:normal-form} follows directly by
contracting a two-vertex word against the two endpoint anchors.  Its action is
\[
 \left(\prod_{j<i}\sigma^j(a_j)\right)
 S^i(x)
 \left(\prod_{j>i}\sigma^j(a_j)\right)c,
\]
with the relation $S^m=L_cR_{c^{-1}}$ used when $i=m$.

Every context is a graft product of elementary contexts.  If normal forms
$L_aR_bS^r$ and $L_dR_eS^s$ are grafted, put
$\delta=\lfloor(r+s)/m\rfloor$.  Relations
\eqref{eq:crossed-relations} give
\begin{equation}\label{eq:normal-multiplication}
 (L_aR_bS^r)(L_dR_eS^s)
 =L_{a\sigma^r(d)c^\delta}
  R_{c^{-\delta}\sigma^r(e)b}
  S^{r+s-\delta m}.
\end{equation}
The same carry $c^\delta$ is obtained from the flattened word because
$\sigma^m=\Ad_c$.  Induction on the number of vertices proves
\eqref{eq:normal-form}.
\end{proof}

\begin{theorem}[Cyclic crossed-product envelope]\label{thm:skew-envelope}
There is a natural isomorphism of $k$-algebras
\begin{equation}\label{eq:skew-envelope}
 \Uac(A,\sigma,n)\cong\mathcal C(A,\sigma,c).
\end{equation}
Consequently endpoint-anchored slot modules are naturally $A$-bimodules $M$
equipped with an invertible $k$-linear operator $S$ such that
\begin{equation}\label{eq:equivariant-bimodule}
 S(axb)=\sigma(a)S(x)\sigma(b),
 \qquad
 S^m(x)=cxc^{-1}.
\end{equation}
For such a module, the $i$-th slot action is
\begin{equation}\label{eq:slot-formula}
 \lambda_i(a_0,\ldots,a_{i-1},x,a_{i+1},\ldots,a_m)
 =p_iS^i(x)q_i,
\end{equation}
where
\[
 p_i=\prod_{j<i}\sigma^j(a_j),
 \qquad
 q_i=\left(\prod_{j>i}\sigma^j(a_j)\right)c.
\]
\end{theorem}

\begin{proof}
Send a context $C$ to
$(a_C\otimes b_C^{\op})u^{r_C}$.  Theorem~\ref{thm:total-assoc} and the
endpoint anchors show that this is well defined, while
\eqref{eq:normal-multiplication} is precisely multiplication in
$\mathcal C(A,\sigma,c)$.  Conversely, send
$(a\otimes b^{\op})u^r$ to $L_aR_bS^r$.  Relations
\eqref{eq:crossed-relations} verify the defining crossed-product relations.
The two maps are inverse by Lemma~\ref{lem:normal-form}.  A module over the
crossed product restricts to an $A^e$-module, hence to an $A$-bimodule, while
$u$ acts by an invertible operator $S$ satisfying
\eqref{eq:equivariant-bimodule}.  Conversely, such a pair $(M,S)$ satisfies
the defining crossed-product relations and therefore extends uniquely to a
$\mathcal C(A,\sigma,c)$-module.  Formula \eqref{eq:slot-formula} is obtained
by evaluating the elementary context with the hole in position $i$.
The long-word formula \eqref{eq:long-word} and $S^m=c(\mathord{-})c^{-1}$
show that all parenthesized one-hole words agree, while
\eqref{eq:c-anchor} follows by direct substitution.  Hence this construction
recovers exactly the endpoint-anchored slot modules.
\end{proof}

\begin{corollary}[Finite-order specialization]\label{cor:finite-order-envelope}
If $c=1$, then $\sigma^m=1$ and
\begin{equation}\label{eq:finite-order-envelope}
 \Ua(A,\sigma,n):=\Uac(A,\sigma,n)
 \cong A^e\rtimes_{\sigma^e}C_m.
\end{equation}
The regular module is anchored, and its remaining unit contexts act as the
powers of $\sigma$.
\end{corollary}

In the finite-order case we abbreviate $T(A,\sigma):=T(A,\sigma,1)$.

\begin{remark}\label{rem:generality}
Theorem~\ref{thm:skew-envelope} separates two sources of higher-arity data.
The quotient group $C_{n-1}$ records the position of the module variable,
while the relation $u^{n-1}=c\otimes(c^{-1})^{\op}$ records the inner defect of
the automorphism period.  The skew-group envelope used below is the untwisted
case $c=1$; the crossed-product theorem remains available for
noncommutative examples with genuinely inner period.
\end{remark}

\begin{example}[A nontrivial central carry]\label{ex:central-carry}
Let $A=k[t,t^{-1}]$, let $\sigma=\id$, and take $c=t$.  Then
\[
 [a_0,\ldots,a_m]=t\,a_0\cdots a_m
\]
is totally associative and endpoint-anchored at $t^{-1}$.  Its envelope is
\[
 (A\otimes_kA) [u]\big/(u^m-t\otimes t^{-1}).
\]
In $A\otimes_kA\cong k[t^{\pm1},s^{\pm1}]$, the monomial $ts^{-1}$ is not an
$m$-th power for $m>1$.  Hence the carry relation cannot be removed by a coefficient rescaling
$u\mapsto vu$.  This gives a concrete crossed-product envelope that is not
reduced to the untwisted presentation by the standard change of generator.
\end{example}

\section{Permutation algebras and orbit--stabilizer blocks}
\label{sec:permutation}

From now on we specialize to $c=1$, so $\sigma^m=1$ and the envelope is the skew group algebra of Corollary~\ref{cor:finite-order-envelope}.

Let $X$ be a finite nonempty $C_m$-set and let $A=k^X$ be the algebra of functions
$X\to k$ with pointwise multiplication.  Write $g$ for a generator of $C_m$
and let $\sigma$ be the induced automorphism.  Then
$A^e\cong k^{X\times X}$, and $C_m$ acts diagonally on $X\times X$.

\begin{theorem}[Orbit--stabilizer envelope]\label{thm:orbit}
For every orbit $O\subseteq X\times X$, choose $x_O\in O$ and put
$H_O=\Stab_{C_m}(x_O)$.  Then, after choosing orbit representatives, there is
an isomorphism
\begin{equation}\label{eq:orbit-decomp}
 \Ua(k^X,\sigma,n)
 \cong
 \prod_{O\in(X\times X)/C_m}
 M_{|O|}(kH_O).
\end{equation}
In particular, the Morita class of the anchored envelope is determined by the
stabilizers of the diagonal orbits.
\end{theorem}

\begin{proof}
By Theorem~\ref{thm:skew-envelope}, the envelope is
$k^{X\times X}\rtimes C_m$, the algebra of the finite action groupoid
$(X\times X)\rtimes C_m$.  The central idempotent
$e_O=\sum_{x\in O}\eps_x$ splits it as a product over orbits.  The action
groupoid on a transitive orbit $O$ is a connected finite groupoid with
isotropy group $H_O$.  Choose, for each $y\in O$, an element $t_y\in C_m$
sending $x_O$ to $y$.  The elements corresponding to arrows
$t_yht_z^{-1}:z\to y$, with $h\in H_O$, multiply as matrix units with
coefficients in $kH_O$.  This identifies the orbit algebra with
$M_{|O|}(kH_O)$.  Taking the product over $O$ proves
\eqref{eq:orbit-decomp}.
\end{proof}

\begin{corollary}[Split orbit blocks]\label{cor:split-orbit-blocks}
Assume that $k$ is a field, $\operatorname{char}k\nmid m$, and $k$ contains all
$m$-th roots of unity.  For each diagonal orbit $O$, let
$H_O=\Stab_{C_m}(x_O)$ and let $\widehat H_O$ be its character group.  Then
\begin{equation}\label{eq:split-orbit-blocks}
 \Ua(k^X,\sigma,n)
 \cong
 \prod_{O\in(X\times X)/C_m}\ \prod_{\chi\in\widehat H_O}
 M_{|O|}(k).
\end{equation}
Consequently the finite-dimensional anchored module category is semisimple
with $\sum_O|H_O|$ simple objects; the simple indexed by $(O,\chi)$ has underlying
$k$-dimension $|O|$.
\end{corollary}

\begin{proof}
Every subgroup $H_O$ is cyclic.  Under the hypotheses on $k$, its group algebra
splits as $kH_O\cong\prod_{\chi\in\widehat H_O}k$.  Apply this to
\eqref{eq:orbit-decomp}.  Morita equivalence for each matrix factor shows that
its unique simple module has dimension $|O|$.
\end{proof}

\begin{proposition}[Simple anchored modules in arbitrary characteristic]
\label{prop:modular-simples}
Let $k$ be algebraically closed of characteristic $p\geq0$.  For a finite
cyclic group $H$, write $|H|_{p'}$ for the order of its maximal prime-to-$p$
quotient, with $|H|_{0'}=|H|$.  Then the finite-dimensional anchored module
category for $T(k^X,\sigma)$ has
\begin{equation}\label{eq:NX}
 N_X:=\sum_{O\in(X\times X)/C_m}|H_O|_{p'}
\end{equation}
isomorphism classes of simple objects.  The simples belonging to an orbit $O$
have underlying $k$-dimension $|O|$.
\end{proposition}

\begin{proof}
The simple $kH_O$-modules are the one-dimensional characters of the maximal
prime-to-$p$ quotient of $H_O$.  Their number is $|H_O|_{p'}$.  Morita equivalence for $M_{|O|}(kH_O)$ multiplies their underlying
vector-space dimension by $|O|$.  Apply Theorem~\ref{thm:orbit} orbit by
orbit.
\end{proof}

\begin{proposition}[Orbit arithmetic]\label{prop:orbit-arithmetic}
Write the orbit decomposition of $X$ as
\[
 X=\coprod_{\alpha\in I}C_m/H_\alpha,
 \qquad \ell_\alpha=[C_m:H_\alpha].
\]
For each ordered pair $(\alpha,\beta)$, the diagonal action on
$C_m/H_\alpha\times C_m/H_\beta$ has
$\gcd(\ell_\alpha,\ell_\beta)$ orbits.  Every such orbit has length
$\operatorname{lcm}(\ell_\alpha,\ell_\beta)$ and stabilizer
$H_\alpha\cap H_\beta$.  Consequently, the gcd $d_X$ of the diagonal
orbit lengths satisfies
\begin{equation}\label{eq:dx-original-action}
 d_X=\gcd_{\alpha\in I}\ell_\alpha
     =[C_m:\langle H_\alpha:\alpha\in I\rangle].
\end{equation}
If $k$ is algebraically closed of characteristic $p\geq0$, then
\begin{equation}\label{eq:NX-original-action}
 N_X=\sum_{\alpha,\beta\in I}
 \gcd(\ell_\alpha,\ell_\beta)\,|H_\alpha\cap H_\beta|_{p'}.
\end{equation}
\end{proposition}

\begin{proof}
Because $C_m$ is abelian, diagonal orbits on
$C_m/H_\alpha\times C_m/H_\beta$ are indexed by the double cosets
$H_\alpha\backslash C_m/H_\beta$, hence by the cosets of
$H_\alpha H_\beta$.  Their number is
$[C_m:H_\alpha H_\beta]=\gcd(\ell_\alpha,\ell_\beta)$.  The stabilizer of
any point is $H_\alpha\cap H_\beta$, so orbit--stabilizer gives orbit length
$[C_m:H_\alpha\cap H_\beta]=\operatorname{lcm}(\ell_\alpha,\ell_\beta)$.
Taking the gcd over all pairs yields
\[
 \gcd_{\alpha,\beta}\operatorname{lcm}(\ell_\alpha,\ell_\beta)
 =\gcd_\alpha\ell_\alpha,
\]
because the diagonal pairs $(\alpha,\alpha)$ occur and every displayed lcm
is divisible by $\gcd_\alpha\ell_\alpha$.  For subgroups of a cyclic group,
the index of the subgroup generated by the $H_\alpha$ is the gcd of their
indices, proving \eqref{eq:dx-original-action}.  Formula
\eqref{eq:NX-original-action} follows from
Proposition~\ref{prop:modular-simples}, since each of the indicated diagonal
orbits has stabilizer $H_\alpha\cap H_\beta$.
\end{proof}

We now specialize to the regular cyclic action.  Identify $X=C_m$ additively,
let $\eps_t$ be the characteristic function of $t$, and set
$\sigma(\eps_t)=\eps_{t+1}$.

\begin{definition}\label{def:regular-cyclic}
The \emph{regular cyclic $n$-ary algebra} $\mathcal T_m(k)$ is
$T(k^{C_m},\sigma)$, where $n=m+1$.
\end{definition}

On primitive idempotents its multiplication is particularly sparse:
\begin{equation}\label{eq:idempotent-product}
 [\eps_{i_0},\ldots,\eps_{i_m}]_\sigma
 =\begin{cases}
   \eps_c,& i_j+j\equiv c\pmod m\text{ for every }j,\\
   0,&\text{otherwise}.
  \end{cases}
\end{equation}

\begin{proposition}[No polyadic unit]\label{prop:no-unit}
If $k$ is nonzero and $m>1$, the algebra $\mathcal T_m(k)$ has no polyadic unit.  In particular,
it is not obtained as the iterated $n$-fold multiplication of a unital binary
algebra on the same underlying vector space.
\end{proposition}

\begin{proof}
Suppose $e\in k^{C_m}$ were a polyadic unit.  Unit action in position $0$
would give
\[
 f\prod_{j=1}^{m}\sigma^j(e)=f
 \qquad(f\in k^{C_m}),
\]
so $P=\prod_{j=1}^{m}\sigma^j(e)$ is the constant function $1$.  Hence every
coordinate of $e$ is a unit.  Unit action in position $1$ would then have the
form
\[
 d\,\sigma(f)=f
 \qquad(f\in k^{C_m})
\]
for an invertible diagonal function $d$.  Taking $f=\eps_t$, the left side is
supported at $t+1$, while the right side is supported at $t$.  This is
impossible for $m>1$.  An iterated product of a unital binary algebra has its
binary unit as a polyadic unit, proving the last assertion.
\end{proof}

The diagonal action on $C_m\times C_m$ has the $m$ free orbits
\[
 O_d=\{(t,t+d):t\in C_m\},
 \qquad d\in C_m.
\]

\begin{corollary}[Matrix-block envelope]\label{cor:matrix-block}
For every commutative ring $k$,
\begin{equation}\label{eq:matrix-block}
 \Ua(\mathcal T_m(k))
 \cong\prod_{d\in C_m}M_m(k).
\end{equation}
If $k$ is a field, the category of finite-dimensional anchored slot modules is
semisimple with exactly $m$ simple objects, each of underlying $k$-dimension
$m$.
\end{corollary}

\begin{proof}
Every orbit $O_d$ has size $m$ and trivial stabilizer, so
Theorem~\ref{thm:orbit} gives \eqref{eq:matrix-block}.  Over a field, each
matrix block has one simple module, namely $k^m$.
\end{proof}

For completeness, matrix units can be written explicitly.  Let
$\eps_{(p,p+d)}$ be the primitive idempotent of $k^{C_m\times C_m}$ and let
$u$ implement the diagonal shift.  Then
\[
 E^{(d)}_{pq}=\eps_{(p,p+d)}u^{p-q}
\]
satisfies
$E^{(d)}_{pq}E^{(e)}_{rs}=\delta_{d,e}\delta_{q,r}E^{(d)}_{ps}$.

\section{Dimension-defect and arity-defect spectra}
\label{sec:K}

We first isolate the mechanism behind the relative calculation.  It applies
beyond the cyclic family and explains why the common dimension of simple
objects produces a Moore-spectrum term.  Throughout, $\K(\mathcal E)$ denotes
connective Quillen $K$-theory of an exact category $\mathcal E$; for a
finite-length module category this is its exact-category, or $G$-theoretic,
$K$-theory rather than the $K$-theory of projective modules; we use the
standard exact-category conventions of \cite{Quillen1973,Buehler2010,Weibel2013}.

\begin{theorem}[Dimension-defect theorem]\label{thm:dimension-defect}
Let $k$ be a field and let $\mathcal A$ be a finite-length $k$-linear abelian
category with simple objects $S_1,\ldots,S_N$ such that
$\End_{\mathcal A}(S_i)=k$.  Let
$F:\mathcal A\to\Vect_k^{\mathrm{fd}}$ be exact and $k$-linear, and put
$w_i=\dim_kF(S_i)>0$ and $d=\gcd(w_1,\ldots,w_N)$.  Then there is a
noncanonical equivalence
\begin{equation}\label{eq:dimension-defect-spectrum}
 \hofib\!\left(\K(\mathcal A)\xrightarrow{\K(F)}\K(k)\right)
 \simeq \K(k)^{N-1}\times\Omega\K(k;\mathbb Z/d).
\end{equation}
In particular,
\begin{equation}\label{eq:dimension-defect-minus-one}
 \pi_{-1}\hofib\K(F)\cong\mathbb Z/d.
\end{equation}
\end{theorem}

\begin{proof}
Quillen d\'evissage identifies $\K(\mathcal A)$ with the $K$-theory of the
semisimple subcategory generated by its simple objects, hence with
$Y^N$, where $Y=\K(k)$.  Under this identification, $\K(F)$ is the weighted
fold map $(w_1,\ldots,w_N):Y^N\to Y$.  Smith normal form gives an integral
change of basis taking this row to $(d,0,\ldots,0)$.  The induced
$\mathrm{GL}_N(\mathbb Z)$-self-equivalence of $Y^N$ yields
\[
 \hofib\K(F)\simeq Y^{N-1}\times\hofib(d:Y\to Y)
 \simeq Y^{N-1}\times\Omega\K(k;\mathbb Z/d).
\]
The splitting depends on the chosen Smith basis and is therefore
noncanonical.  Since $Y$ is connective and $K_0(k)=\mathbb Z$, taking
$\pi_{-1}$ gives \eqref{eq:dimension-defect-minus-one}.
\end{proof}

We now apply the theorem to a finite $C_m$-set $X$.  Let $k$ be algebraically
closed of characteristic $p\geq0$, let $\mathcal M_X(k)$ be the exact category
of finite-dimensional anchored slot modules over $T(k^X,\sigma)$, and let
$F_X:\mathcal M_X(k)\to\Vect_k^{\mathrm{fd}}$ forget the slot action.  Define
\begin{equation}\label{eq:orbit-defect}
 \R_X(k):=\hofib\!\left(
 \K(\mathcal M_X(k))\xrightarrow{\K(F_X)}\K(k)
 \right),
\end{equation}
and put
\begin{equation}\label{eq:dX}
 d_X:=\gcd\{|O|:O\in(X\times X)/C_m\}.
\end{equation}
For a positive integer $d$, write
$\K(k;\mathbb Z/d):=\cofib(d:\K(k)\to\K(k))$.

\begin{corollary}[Orbit-defect spectrum]\label{thm:orbit-defect}
With $N_X$ as in \eqref{eq:NX}, there is a noncanonical equivalence
\begin{equation}\label{eq:orbit-defect-spectrum}
 \R_X(k)\simeq
 \K(k)^{N_X-1}\times\Omega\K(k;\mathbb Z/d_X).
\end{equation}
Moreover,
\begin{equation}\label{eq:orbit-minus-one}
 \pi_{-1}\R_X(k)\cong\mathbb Z/d_X
 \cong
 \mathbb Z/\gcd_\alpha\ell_\alpha,
\end{equation}
where $X=\coprod_\alpha C_m/H_\alpha$ and
$\ell_\alpha=[C_m:H_\alpha]$.  Hence the boundary vanishes exactly when the
point stabilizers generate $C_m$, and it has maximal order $m$ exactly when
the action on $X$ is free.
\end{corollary}

\begin{proof}
By Proposition~\ref{prop:modular-simples}, the simple objects indexed by a
diagonal orbit $O$ have endomorphism ring $k$ and underlying vector-space
dimension $|O|$.  Apply Theorem~\ref{thm:dimension-defect}; the gcd of the
weights is $d_X$.  Proposition~\ref{prop:orbit-arithmetic} gives the final
form of the boundary group.
\end{proof}

\begin{corollary}[Divisor realization invisible to ordinary $K$-theory]
\label{cor:ordinary-invisible}
Let $k$ be algebraically closed of characteristic $0$.  For every divisor
$d\mid m$, let $X_d=C_m/H_d$ be the transitive $C_m$-set of length $d$.
Then
\begin{equation}\label{eq:divisor-family-ordinary}
 \K(\mathcal M_{X_d}(k))\simeq\K(k)^m
 \qquad(d\mid m),
\end{equation}
whereas
\begin{equation}\label{eq:divisor-family-relative}
 \pi_{-1}\R_{X_d}(k)\cong\mathbb Z/d.
\end{equation}
Thus the divisor lattice of $m$ is realized by relative boundary groups even
though all the ordinary exact $K$-theory spectra in the family are equivalent.
In particular, the one-point and regular actions give the extreme defects
$0$ and $\mathbb Z/m$.
\end{corollary}

\begin{proof}
The diagonal action on $X_d\times X_d$ has $d$ orbits, each with stabilizer
$H_d$ of order $m/d$.  In characteristic $0$, each block therefore contributes
$m/d$ simple modules, so $N_{X_d}=d(m/d)=m$.  Proposition~\ref{prop:orbit-arithmetic}
gives $d_{X_d}=d$.  Apply Corollary~\ref{thm:orbit-defect}.
\end{proof}

\begin{corollary}[Stable defect classification]\label{cor:defect-classification}
For algebraically closed $k$, the relative spectrum $\R_X(k)$ determines the
pair $(N_X,d_X)$: one has
\[
 N_X=1+\operatorname{rank}\pi_0\R_X(k),
 \qquad
 d_X=|\pi_{-1}\R_X(k)|,
\]
with the convention that the trivial group has order $1$.  Conversely,
$(N_X,d_X)$ determines $\R_X(k)$ up to the noncanonical equivalence in
\eqref{eq:orbit-defect-spectrum}.
\end{corollary}

\begin{proof}
Multiplication by $d_X$ is surjective on $K_1(k)=k^\times$ because $k$ is
algebraically closed, and it has zero kernel on $K_0(k)=\mathbb Z$.  Hence the
Moore-loop factor in \eqref{eq:orbit-defect-spectrum} contributes neither a
free summand nor torsion to $\pi_0$.  The assertions follow from
\eqref{eq:orbit-defect-spectrum} and \eqref{eq:orbit-minus-one}.
\end{proof}

\begin{remark}\label{rem:orbit-meaning}
By \eqref{eq:dx-original-action}, the boundary measures the residual cyclic
period not generated by point stabilizers.  A fixed point has stabilizer
$C_m$, hence $d_X=1$ and zero negative boundary; a free action has trivial
stabilizers and $d_X=m$.  Thus \eqref{eq:orbit-minus-one} is an invariant of
the positional action, not a formal rereading of the number of inputs.
\end{remark}

\begin{example}[Intermediate defect]\label{ex:intermediate-defect}
Let $k$ be algebraically closed of characteristic $0$, and let $C_6$ act on
\[
 X=C_6/C_3\ \coprod\ C_6.
\]
The two orbit lengths are $2$ and $6$.  Proposition~\ref{prop:orbit-arithmetic}
gives $d_X=2$ and
\[
 N_X=2\cdot3+2+2+6=16.
\]
Hence
\[
 \K(\mathcal M_X(k))\simeq\K(k)^{16},
 \qquad
 \R_X(k)\simeq\K(k)^{15}\times\Omega\K(k;\mathbb Z/2).
\]
Thus the boundary may be a proper nontrivial divisor of $m$; here it records
the residual index of the subgroup $C_3$ generated by the point stabilizers.
\end{example}

We now return to the regular cyclic algebra $\mathcal T_m(k)$.  No restriction
on the characteristic is needed because its anchored envelope is already a
product of matrix algebras.  Let $\mathcal M_m(k)$ be its finite-dimensional
anchored module category, let
$F_m:\mathcal M_m(k)\to\Vect_k^{\mathrm{fd}}$ be forgetful, and define
\begin{equation}\label{eq:arity-defect}
 \R_m(k):=\hofib\!\left(
 \K(\mathcal M_m(k))\xrightarrow{\K(F_m)}\K(k)
 \right).
\end{equation}

\begin{corollary}[Arity detection]\label{thm:arity-detection}
For every field $k$ and $m\geq2$, there is a noncanonical equivalence
\begin{equation}\label{eq:main-spectrum}
 \R_m(k)
 \simeq
 \K(k)^{m-1}\times\Omega\K(k;\mathbb Z/m).
\end{equation}
In particular,
\begin{equation}\label{eq:minus-one}
 \pi_{-1}\R_m(k)\cong\mathbb Z/m.
\end{equation}
Thus the relative spectrum determines the arity of $\mathcal T_m(k)$ via
\[
 n=1+|\pi_{-1}\R_m(k)|.
\]
\end{corollary}

\begin{proof}
By Corollary~\ref{cor:matrix-block}, the module category is the product of $m$
copies of $\Vect_k^{\mathrm{fd}}$, and each standard simple has underlying
dimension $m$.  Therefore $\K(F_m)=m\nabla:\K(k)^m\to\K(k)$.  Equivalently,
apply Theorem~\ref{thm:orbit-defect} formally to the free diagonal orbits, for
which $N_X=d_X=m$; its proof uses only the matrix-block decomposition and so
remains valid in every characteristic.
\end{proof}

\begin{remark}\label{rem:why-relative}
The ordinary exact $K$-theory of $\mathcal M_m(k)$ is $\K(k)^m$ and records the
number of matrix blocks.  The forgetful map also records the common dimension
$m$ of the simple anchored modules.  Its homotopy fibre combines the block
count with this dimension, and the degree $-1$ boundary is the obstruction to
$m\mathbb Z\subset\mathbb Z$ being surjective.  This boundary is what detects
the arity in the regular family.
\end{remark}

\begin{corollary}\label{cor:long-exact}
For every $i\geq0$ there is a natural short exact sequence
\begin{equation}\label{eq:short-exact}
\begin{aligned}
 0&\longrightarrow
 \operatorname{coker}\bigl(m:K_{i+1}(k)\to K_{i+1}(k)\bigr)\\
 &\longrightarrow \pi_i\hofib(m)
 \longrightarrow
 \ker\bigl(m:K_i(k)\to K_i(k)\bigr)
 \longrightarrow0.
\end{aligned}
\end{equation}
and
\[
 \pi_i\R_m(k)\cong K_i(k)^{m-1}\oplus\pi_i\hofib(m).
\]
\end{corollary}

\begin{proof}
Apply homotopy groups to the fibre sequence of multiplication by $m$ and use
\eqref{eq:main-spectrum}.
\end{proof}

\section{Finite fields and degree-wise torsion formulas}
\label{sec:finite-fields}

Let $k=\mathbb F_q$.  Quillen proved
\cite{Quillen1972}
\begin{equation}\label{eq:finite-field-K}
 K_0(\mathbb F_q)=\mathbb Z,
 \qquad
 K_{2j}(\mathbb F_q)=0,
 \qquad
 K_{2j-1}(\mathbb F_q)\cong\mathbb Z/(q^j-1)
 \quad(j\geq1).
\end{equation}
Put $g_j=\gcd(m,q^j-1)$.

\begin{theorem}[Finite-field calculation]\label{thm:finite-field}
For the regular cyclic $(m+1)$-ary algebra $\mathcal T_m(\mathbb F_q)$,
\begin{align}
 \pi_{-1}\R_m(\mathbb F_q)
 &\cong\mathbb Z/m,\label{eq:ff-minus1}\\
 \pi_0\R_m(\mathbb F_q)
 &\cong\mathbb Z^{m-1}\oplus\mathbb Z/g_1,\label{eq:ff-zero}\\
 \pi_{2j-1}\R_m(\mathbb F_q)
 &\cong
 (\mathbb Z/(q^j-1))^{m-1}\oplus\mathbb Z/g_j
 \qquad(j\geq1),\label{eq:ff-odd}\\
 \pi_{2j}\R_m(\mathbb F_q)
 &\cong\mathbb Z/g_{j+1}
 \qquad(j\geq1).\label{eq:ff-even}
\end{align}
\end{theorem}

\begin{proof}
For a cyclic group $\mathbb Z/N$, both the kernel and cokernel of
multiplication by $m$ are cyclic of order $\gcd(m,N)$.  Insert
\eqref{eq:finite-field-K} into Corollary~\ref{cor:long-exact}.  In every
positive degree one of the two adjacent finite-field $K$-groups vanishes, so
the short exact sequence has only one nonzero end term.  Degree $0$ also
contains the factor $K_0(\mathbb F_q)^{m-1}=\mathbb Z^{m-1}$, while degree
$-1$ is Theorem~\ref{thm:arity-detection}.
\end{proof}

\begin{corollary}[Prime-wise periodicity]\label{cor:torsion-criterion}
Let $\ell$ be a prime divisor of $m$.
\begin{enumerate}[label=\textup{(\roman*)}]
\item If $\ell\mid q$, then the additional coefficient terms in
\eqref{eq:ff-zero}--\eqref{eq:ff-even} have trivial $\ell$-primary part.
\item If $\ell\nmid q$ and $o_\ell(q)$ is the multiplicative order of $q$
modulo $\ell$, then the additional $\ell$-primary term in degree $2j-1$ is
nonzero exactly when $o_\ell(q)\mid j$, and the one in degree $2j$ is
nonzero exactly when $o_\ell(q)\mid j+1$.
\item In either nonzero case its order is
\[
 \ell^{\min\{v_\ell(m),v_\ell(q^r-1)\}},
\]
with $r=j$ in odd degree and $r=j+1$ in even degree.
\end{enumerate}
Consequently the extra summands are nonzero precisely when
\[
 \gcd(m,q^r-1)>1.
\]
\end{corollary}

\begin{proof}
The $\ell$-primary part of $\mathbb Z/\gcd(m,q^r-1)$ has the displayed order.
If $\ell\mid q$, then $q^r-1\equiv-1\pmod\ell$.  If $\ell\nmid q$, then
$\ell\mid q^r-1$ exactly when $o_\ell(q)\mid r$.  Apply
Theorem~\ref{thm:finite-field}.
\end{proof}

The first examples are summarized below.
\begin{center}
\begin{tabular}{@{}cccc@{}}
\toprule
Arity $n$ & $m=n-1$ & Anchored envelope & $\pi_{-1}\R_m(k)$\\
\midrule
$3$ & $2$ & $M_2(k)\times M_2(k)$ & $\mathbb Z/2$\\
$4$ & $3$ & $M_3(k)^3$ & $\mathbb Z/3$\\
$5$ & $4$ & $M_4(k)^4$ & $\mathbb Z/4$\\
$6$ & $5$ & $M_5(k)^5$ & $\mathbb Z/5$\\
\bottomrule
\end{tabular}
\end{center}

\begin{example}\label{ex:ternary}
For $n=3$, let $A=k\times k$ and let $\sigma$ exchange the two coordinates.
Then
\[
 [a,b,c]_\sigma=a\,\sigma(b)c
\]
is totally associative and has no ternary unit.  Its anchored envelope is
$M_2(k)\times M_2(k)$, and
\[
 \R_2(k)\simeq\K(k)\times\Omega\K(k;\mathbb Z/2).
\]
Over $\mathbb F_q$, the relative boundary is $\mathbb Z/2$, while
$\pi_0$ has an additional $\mathbb Z/2$ precisely when $q$ is odd.
\end{example}

\begin{example}\label{ex:quaternary}
For $n=4$, take $A=k^3$ with the cyclic shift.  The product
\[
 [a,b,c,d]_\sigma=a\,\sigma(b)\sigma^2(c)d
\]
has anchored envelope $M_3(k)^3$ and
\[
 \R_3(k)\simeq\K(k)^2\times\Omega\K(k;\mathbb Z/3).
\]
Thus $\pi_{-1}\cong\mathbb Z/3$, and over $\mathbb F_q$ the coefficient
summands occur in the degrees for which $3\mid q^j-1$.
\end{example}

\section{Relation with existing enveloping constructions}
\label{sec:comparison}

The results above combine three established theories in a way that produces a
new computable invariant.

First, automorphism-derived $n$-ary products and representations of polyadic
groups precede this work \cite{DudekGlazek2008,BorowiecDudekDuplij2006,
DudekShahryari2012}, and the abstract existence of an associative envelope
representing modules is part of operadic module theory
\cite{Fresse2009,LodayVallette2012,Khoroshkin2022}.  Our contribution is not the product or abstract representability.  It is the
explicit carried context normal form and the identification of the
endpoint-anchored envelope with the cyclic crossed product
$A^e\langle u\rangle/(ux-\sigma^e(x)u,\,u^{n-1}-z_c)$.  The skew group
algebra $A^e\rtimes C_{n-1}$ is the finite-order specialization $c=1$.  This identification makes
available the structure theory of group actions on associative algebras and
leads directly to the orbit formula \eqref{eq:orbit-decomp}.

Second, Sitarz associates a $\mathbb Z/(n-1)$-graded universal cover to an
$n$-ary associative algebra \cite{Sitarz1998}.  That construction is designed
to contain the original $n$-ary algebra as a homogeneous subspace.  The
anchored envelope is instead characterized by its module category.  For the
inner-periodic cyclic twist, its degree operator becomes the crossed-product
generator $u$; the relations $ux=\sigma^e(x)u$ and $u^m=z_c$ record,
respectively, covariance with the two-sided action and the inner period.  In
the finite-order case this is the group element of a skew group algebra.

Third, cyclic crossed products and skew group algebras have a substantial
representation theory \cite{ReitenRiedtmann1985,Montgomery1993}.  The
orbit--stabilizer theorem is
the finite action-groupoid decomposition specialized to the diagonal action
arising from slot positions.  The regular cyclic action is free, which is why
all blocks are full matrix algebras even when the characteristic divides
$m$.

Finally, the spectra in \eqref{eq:orbit-defect-spectrum} and
\eqref{eq:main-spectrum} are not relative $K$-theory spectra of ideals, nor are
they the polyadic Grothendieck groups of \cite{Duplij2022}.  They are homotopy
fibres of exact forgetful functors.  They also differ from the formal
$\Gamma$-semiring spectra in
\cite{Gokavarapu2025Higher,Gokavarapu2025Fundamental}: the present category is
identified with modules over an explicit cyclic crossed product (and, for
permutation algebras, an explicit skew group algebra), while the forgetful
map itself is computed.  Theorem~\ref{thm:dimension-defect} isolates the Smith-normal-form
mechanism for any finite-length category with a dimension-valued exact
functor.  For a permutation algebra,
Proposition~\ref{prop:orbit-arithmetic} computes the boundary directly from
the orbit lengths of the original cyclic action; Corollary~\ref{cor:ordinary-invisible} shows that this refines
ordinary exact $K$-theory.  For the regular cyclic family the gcd is $n-1$, so
arity becomes a stable homotopy obstruction.

The construction suggests several further directions.  For a nonfree cyclic
action on $X$, Theorem~\ref{thm:orbit} replaces matrix blocks over $k$ by
matrix blocks over stabilizer group algebras; modular representation theory
should then produce extension gluing and singular $K$-theory.  For a noncommutative algebra $A$, the cyclic crossed product in Theorem~\ref{thm:skew-envelope} connects endpoint-anchored slot modules with equivariant Hochschild theory.  It is also natural to ask
for analogous formulas when the cyclic automorphism is replaced by a periodic
anti-automorphism or by a cocycle-twisted action.

\section{Conclusion}

The inner-periodic product \eqref{eq:twisted} provides a totally associative
higher-arity structure in which both the modulus $n-1$ and the inner period
of the automorphism survive passage to modules.  Endpoint anchoring leaves a
cyclic operator whose covariance with left and right multiplication produces
the crossed-product envelope
\[
 A^e\langle u\rangle/
 (ux-\sigma^e(x)u,\ u^{n-1}-c\otimes(c^{-1})^{\op}).
\]
For $c=1$ this specializes to $(A\otimes A^{\op})\rtimes C_{n-1}$.
This turns positional $n$-ary representation theory into equivariant bimodule
theory.

For permutation algebras, the envelope is controlled by diagonal orbits and
stabilizers.  The regular cyclic family is unitless in the polyadic sense and
has envelope $M_{n-1}(k)^{n-1}$.  The corresponding relative exact
$K$-theory is not arity-blind: the Moore-coefficient factor gives
\[
 \pi_{-1}\cong\mathbb Z/(n-1),
\]
while Quillen's finite-field calculation yields degree-by-degree formulas
and prime-wise periodicity criteria in all nonnegative degrees.  More generally,
transitive actions of every divisor length $d\mid n-1$ have the same ordinary
exact $K$-theory but relative boundary $\mathbb Z/d$.  Thus the higher-arity
operation, its representation category, and its relative $K$-theory are linked
by a single cyclic symmetry.

\section*{Declarations}

\noindent\textbf{Funding.} No external funding was received for this work.

\noindent\textbf{Conflict of interest.} The authors declare that they have no conflict of interest.

\noindent\textbf{Data availability.} No datasets were generated or analysed in this theoretical study.

\noindent\textbf{Author contributions.} Chandrasekhar Gokavarapu: conceptualization, initial formulation, manuscript preparation, and coordination. Sajani Lavanya Madasi: algebraic verification, examples, and manuscript review. Madhusudhana Rao Dasari: supervision, mathematical validation, and critical revision. All authors accept responsibility for their stated contributions and for the submitted version.

\end{document}